\documentclass[12pt]{article}

\usepackage{amsmath}
\usepackage{amsfonts}
\usepackage{amssymb}

\newtheorem{proposition}{Proposition}[section]
\newtheorem{theorem}{Theorem}[section]

\newtheorem{remark}[theorem]{Remark}

\def\phi{{\varphi}}

\DeclareSymbolFont{AMSb}{U}{msb}{m}{n}
\DeclareMathSymbol{\N}{\mathbin}{AMSb}{"4E}
\DeclareMathSymbol{\Z}{\mathbin}{AMSb}{"5A}
\DeclareMathSymbol{\R}{\mathbin}{AMSb}{"52}
\DeclareMathSymbol{\Q}{\mathbin}{AMSb}{"51}
\DeclareMathSymbol{\I}{\mathbin}{AMSb}{"49}
\DeclareMathSymbol{\C}{\mathbin}{AMSb}{"43}

\begin{document}
\title{Photo-acoustic tomography in a rotating measurement setting}

\author{{Guillaume Bal\footnote{Department of Applied Physics and Applied Mathematics, Columbia University, New York, NY, 10027, USA. E-mail: gb2030@columbia.edu.}
\qquad Amir Moradifam \footnote{Department of Mathematics, University of California, Riverside, CA, USA. E-mail: moradifam@math.ucr.edu. }}}
\date{\today}
\smallbreak \maketitle

\abstract{Photo-acoustic tomography (PAT) aims to leverage the photo-acoustic coupling between optical absorption of light sources and ultrasound (US) emission to obtain high contrast reconstructions of optical parameters with the high resolution of sonic waves. Quantitative PAT often involves a two-step procedure: first the map of sonic emission is reconstructed from US boundary measurements;  and second optical properties of biological tissues are evaluated. We consider here a practical measurement setting in which such a separation does not apply. We assume that the optical source and an array of ultrasonic transducers are mounted on a rotating frame (in two or three dimensions) so that the light source rotates at the same time as the US measurements are acquired. As a consequence, we no longer have the option to reconstruct a map of sonic emission corresponding to a given optical illumination. We propose here a framework where the two steps are combined into one and an absorption map is directly reconstructed from the available ultrasound measurements.  }

%\title{Photoacoustic tomography with rotating partial measurements}
%
%\author{{Guillaume Bal\footnote{Department of Applied Physics and Applied Mathematics, Columbia University, New York, NY, 10027, USA. E-mail: gb2030@columbia.edu.}
%\qquad Amir Moradifam \footnote{Department of Mathematics, University of California, Riverside, CA, USA. E-mail: moradifam@math.ucr.edu. }}}
%\date{\today}
%\smallbreak \maketitle
%
%\abstract{To be written in later versions. }

\section{Introduction}

Photo-acoustic tomography (PAT) is a novel medical imaging modality that aims to image the optical properties of biological tissues with high resolution. It combines the high contrast of  optical (mostly absorption) parameters with the high resolution of ultrasound. As radiation propagates into tissues, a small fraction is transformed into sonic waves by the photo-acoustic effect. These sonic waves propagate  through the domain and are measured by an array of transducers at the boundary of the domain of interest $\Omega\subset\R^n$, where $n$ is spatial dimension.

Mathematically, sound propagation is modeled by the following scalar wave equation:
\begin{eqnarray}
\left\{ \begin{array}{ll}
(\partial^2_t-c^2\Delta)v=0 &\text{in } (0,\infty) \times \R^n\\
v|_{t=0}=H(x) &\text{in } \R^n\\
\partial_t v|_{t=0}=0 &\text{in } \R^n.
\end{array} \right.
\end{eqnarray}
Here the function $c(x)$ is assumed to be a smooth function and $H(x)$ is the amount of acoustic signal generated by the absorbed radiation of a short pulse of light propagating throughout the domain. Its expression is given by 
\begin{equation}
H(x)=\lambda(x)\sigma(x)u(x),
\end{equation}
where for each $x\in\Omega$, $u(x)$ is the density of  radiation reaching point $x$, $\sigma(x)$ is the absorption coefficient, and $\lambda(x)$ is the Gr\"{u}neisen coefficient, which characterizes the amount of ultrasound generated by each absorbed photon.

A reasonable model for the propagation of radiation is given by the following second-order elliptic equation
\begin{eqnarray}\label{pde1}
\left\{ \begin{array}{ll}
-\nabla \cdot \gamma(x) \nabla u+\sigma u=0 &\text{in } \Omega\\
u=f&\text{on }\partial \Omega, 
\end{array} \right.
\end{eqnarray}
where $\Omega$ is a smooth bounded domain in $\R^n$ ($n\geq 2$),  $\gamma(x)$ is a (scalar) diffusion coefficient and $f$ describes how light enters the domain $\Omega$.

Photoacoustic Tomography aims to reconstruct $(\gamma(x),\sigma(x))$ as well as possibly $\lambda(x)$  from the measurements $v(t,x)|_{\partial \Omega}$ of the acoustic pressure leaving the domain $\Omega$ at each point $x\in\partial \Omega$ and each positive time $t>0$; as well as for each available probing illumination $f$.

\medskip

Such reconstructions are typically done in two steps. In a first step, the sonic source $H$ is reconstructed from the measurements $v(t,x)|_{\partial \Omega}$. This is a well-posed problem, at least in the presence of complete data and when sound speed $c(x)$ is known; see, e.g., \cite{LU, JSUZ, SU13, SU09,W-JDM-04,XW-RSI-06} and their references for experimental and theoretical works on the first, qualitative, step of PAT. 

Once $H(x)$ has been reconstructed for one or more illuminations $f$ on $\partial\Omega$, the second, quantitative, step of PAT allows us to obtain explicit reconstructions of $\gamma$, $\sigma$, and $\lambda$ in some cases. We refer to \cite{BR11, BU13,CAB-JOSA-09,CLB-SPIE-09,CLB-BOE-10,Z-AO-10} as well as their references for several works on the problem. Note that high resolution reconstructions can typically not be achieved by purely optical measurements, which are modeled by a problem that has similar stability properties to the standard electrical impedance tomography problem; see \cite{AS-IP-10}.

\medskip

In this paper, we consider an experimental setting  \cite{BSFECO-JBO-09,LNWEOA-SPIE-15} in which such a separation into first and second steps is not feasible. The reason is that the source of radiation $f$ and the small array of transducers performing the acoustic measurements are mounted on a rotating frame. In other words, as the pressure $v(t,x)$ is measured by a rotating array of transducers on $\partial\Omega$, it corresponds to a radiation source $f$ that also rotates.  We thus no longer acquire a pressure $v(t,x)$ that corresponds to a single illumination $f$ and cannot even define a meaningful initial condition $H(x)$. The objective of this paper is to present a mathematical framework for this experimental setting, whose main advantage is that it allows for a clear spatial separation between the light source and the array of detectors. 

As we mentioned above, the qualitative and quantitative steps of PAT need to be merged into one reconstruction. In order to simplify the presentation, we assume that the Gr\"uneisen coefficient is known and set to $1$. We also assume that the diffusion coefficient $\gamma$ is known and is normalized to $1$. Under these assumptions, we present a theory for the reconstruction of the absorption coefficient $\sigma$ from knowledge of pressure measurements at the domain's boundary in the aforementioned rotating setting.

The rest of the paper is structured as follows. The measurement setting and our main results are presented in Section \ref{sec:main}. The proof of our main result is split into a proof of a linearized version in Section \ref{sec:lin} and a proof of the full nonlinear inverse problem in Section \ref{sec:nlin}.

\section{Measurement setting and main results}
\label{sec:main}
%Photoacoustic tomography (PAT) is a hybrid imaging method that combines high contrast of optical waves with the high resolution of acoustic waves in order to obtain images with high resolution and high contrast. In PAT a tissue is exposed to a short pulse radiation and absorbs a small fraction of the radiation which leads to a slight temperature increase resulting in a mechanical expansion. Such expansion is sufficient to emit acoustic waves that travel back to the boundary of the tissue where they are measured by a set of transducers. 
%
%
%Let $u(x)$ be the total intensity of radiation reaching a point $x$ over the time span of radiation, then $u$ satisfies 
%
%\begin{eqnarray}\label{pde1}
%\left\{ \begin{array}{ll}
%-\nabla \cdot \gamma(x) \nabla u+\sigma u=0 &\text{in } \Omega\\
%u=f&\text{on }\partial \Omega, 
%\end{array} \right.
%\end{eqnarray}
%where $\Omega$ is a smooth bounded domain in $\R^n$ ($n\geq 2$), $\gamma(x)$ is a scalar diffusion coefficient and $\sigma(x)$ is the absorption coefficient. The amount of energy transformed into acoustic waves is given by 
%\begin{equation}
%H(x)=\lambda(x)\sigma(x)u(x),
%\end{equation}
%where $\lambda(x)$ is the Gr\"{u}neisen coefficient and $\sigma(x)u(x)$ is the density of absorbed radiation. The resulting acoustic waves are modeled by 
%
%\begin{eqnarray}
%\left\{ \begin{array}{ll}
%(\partial^2_t-c^2\Delta)v=0 &\text{in } (0,\infty) \times \R^n\\
%v|_{t=0}=H(x) &\text{in } \R^n\\
%\partial_t v|_{t=0}=0 &\text{in } \R^n.
%\end{array} \right.
%\end{eqnarray}

%
%
%
%

We consider the stable reconstruction of the absorption coefficient $\sigma(x)$ from multiple partial measurements under the assumption $\gamma =\lambda= 1$ inside a ball in $\R^n$ ($n\geq 2$). More precisely, let $B_\rho$ be the ball of radius $\rho$ in $\R^n$ and assume $\Omega \Subset B_\rho$. Let $\sigma \in W^{1,\infty}_0(\Omega)$, and define   
\[\Theta:=\{R_i\in SO(n), \ \ i=1,2,...,m\},\]
a finite number of rotations around the origin in $\R^n$. Now fix $0\leq g\in H^{1/2}(\partial B_\rho)\cap L^{\infty}(\partial B_\rho)$ and for $R_i \in \Theta$, define $u_{i}(x)$ to be the unique solution of   
\begin{eqnarray}\label{pde1}
\left\{ \begin{array}{ll}
-\Delta u_i+\sigma u_i=0 &\text{in } B_\rho\\
u=g_{i}&\text{on }\partial B_\rho, 
\end{array} \right.
\end{eqnarray}
where $g_i(x) =g(R_i x)$, $R_i \in \Theta$. We shall assume that $g$ is not identically zero and non-negative. 
In practical applications, we may envision $g$ to have a small support on $\partial B_\rho$, and certainly to be supported away from the location of the ultrasound transducers we now consider. 

Let $v_{i}$ be the solution of the wave equation 
\begin{eqnarray}\label{waveEq}
\left\{ \begin{array}{ll}
(\partial^2_t-c^2\Delta)v=0 &\text{in } (0,\infty) \times \R^n\\
v|_{t=0}=\sigma u_{i} \\
\partial_t v|_{t=0}=0,
\end{array} \right.
\end{eqnarray}
where $c\geq c_0>0$ is the  sound speed, and $c-1$ is assumed to be supported in $\bar{B}_{\rho}$. In general the
time-dependent wave solution at the boundary of the ball $B_\rho$ is given by
\begin{equation}\label{measurmentOperators}
 \Lambda_i^{*} \sigma:=\Lambda(\sigma u_i) :=  v_{i}|_{(0,\infty)\times \partial B_\rho}, \ \quad 1\leq i\leq m. %i\in \Theta.
\end{equation}

In practical setting that we consider here, we have access to $ v_{i}$ on $[0,T] \times \Gamma$ for some $\Gamma \Subset \partial B_\rho$ and not on the whole domain $(0,\infty)\times \partial B_\rho$. Typically the support of the ultrasound transducers $\Gamma$ is relatively small and away from the support of the optical source to avoid measurement interferences.  

To model this restriction, fix $\Gamma \subset \partial B_\rho$ and define $\Gamma_i=R_i(\Gamma)$ for $R_i\in \Theta$. Here, we are interested to know if the absorption coefficient $\sigma(x)$ can be stably determined from the finite number of rotating partial measurements 
\begin{equation}\label{measurements}
\chi_{i}  \Lambda_{i}^*\sigma =\chi_{i}  \Lambda_{i} u_i\sigma= \chi_{i} v_{i}|_{(0,\infty)\times \partial B_{\rho}}, \ \ R_i \in \Theta,
\end{equation}
where $\chi_{i} \in C^{\infty}_0([0,\infty)\times B_\rho)$ are cut-off functions with supp$(\chi_i) \subset [0,\infty)\times \Gamma_i$. This is a reasonably faithful model for the experimental setups described in \cite{BSFECO-JBO-09,LNWEOA-SPIE-15}, to which we refer the reader for additional details.

\medskip

We now introduce additional hypotheses and notation in order to state our main result; theorem \ref{MainTheorem} below. One of our main theoretical tools is a description of acoustic wave propagation in the domain $B_\rho$ mostly following the presentation in \cite{SU09}.
As in \cite{SU09}, and for a given $h$, we define $v$ to be the unique solution of 
\begin{eqnarray}
\left\{ \begin{array}{cl}
(\partial^2_t-c^2\Delta)v &= 0 \ \ \text{in } (0,T) \times \R^n\\
v|_{[0,T]\times \partial B_{\rho}} &= h,\\
v|_{t=T} &=\varphi,\\
\partial_t v|_{t=T } &=0,
\end{array} \right.
\end{eqnarray}
where $\varphi$ is the harmonic extension of $h(T, \cdot )$ in $B_\rho$. We also define
\begin{equation}\label{A}
Ah:= v(0, \cdot) \ \ \hbox{in} \ \ \bar{B_\rho},
\end{equation}
the wave solution at time $t=0$, and set 
\begin{equation}
\mathcal{G}:=\{(t,y): \ \ y\in \cup_{R_i \in \Theta}\Gamma_i, 0<t<s(y)\},
\end{equation}
where $s(y)$ is a continuous function on $\cup_{R_i \in \Theta}\Gamma_i$ indicating how long measurements need to last at every measurement point. It is know from, e.g., the work in \cite{SU09} that longer times are necessary for stability purposes than for injectivity purposes.

As in (\cite{SU09}) and to guarantee injectivity of the measurement operator, we will assume that there exists $j \in \{1,2,...,m\}$ such that 
\begin{equation}\label{Uniqueness.con}
\forall x\in \Omega, \ \ \exists y\in \Gamma_j \ \ \hbox{with} \ \ \hbox{dist}(x,y)<s(y), 
\end{equation}
where dist$(x,y)$ denotes the distance with respect to the metric $c^{-2}g$ ($g$ is the Euclidean metric in $\R^n$). 
This assumption is partially technical: it imposes that the measurements for one of the rotation step $j$ are taken for a sufficiently long duration that the measurement operator (mapping $\sigma$ to the available measurements) is injective. 

To guarantee stability of the reconstruction of $\sigma$, we  need the following stronger assumption on $s(y)$:
\begin{equation}\label{stab.con}
\forall(x,\xi) \in WF(\sigma) \cap (\Omega \times \R^n), \ \ (\tau_{\eta}(x,\xi), \gamma_{x,\xi}(\tau_{\eta}(x,\xi)))\in \mathcal{G} \hbox{ for} \ \ \eta=+ \ \ \hbox{or}\ \ \eta=-,
\end{equation}
where $\gamma_{\eta}(x,\xi)$ are the integral curves of the corresponding Hamilton vector field (see Chapter 6 in \cite{Grigis}) and 
\[\tau_{\pm}(x,\xi)=\max \{t \geq 0: \ \ \gamma_{x,\xi}(\pm t)\in \bar{B}_\rho\}.\] 
We refer to \cite{SU09} for additional details on this assumption, which here simply means that any singularity of $\sigma$ at position $x\in B_\rho$ and in direction $\xi$ propagates to a singularity in the measurement set ${\cal G}$.

Note that since $\overline{\Omega}$ is compact there exists an open set $\mathcal{G'} \subset \mathcal{G}$ such that (\ref{stab.con}) still holds. Define
\[\Gamma'_{i}:=\Gamma_{i} \cap \pi_{x} (\mathcal{G'})\]
where $\pi_{x}$ is the projection map. Now fix the cut-off functions $\chi_{i} \in  C^{\infty}_0([0,\infty)\times B_\rho)$ such that 
\[\hbox{supp} (\chi_{i}) \subset \Gamma_{i}\times (0,\infty) \cap \mathcal{G} \ \  \hbox{and} \ \ \chi_{i}=1 \ \ \hbox{on}  \ \ \Gamma'_{i}\times (0,\infty) \cap \mathcal{G}.  \]
The following is the main result of this paper.

\begin{theorem}\label{MainTheorem}
Let $B_{\rho}$ be the ball of radius $\rho$ in $\R^n$, $\Omega \Subset B_\rho$, and $\sigma \in H^{1}_0(\Omega)$. Let $m\in \N$ and assume that (\ref{Uniqueness.con}) and (\ref{stab.con}) hold. If

\[\Lambda_i^* (\sigma)=\Lambda_i^*(\tilde{\sigma}) \ \ \hbox{on} \ \ \Gamma_i\times (0,\infty) \cap\mathcal{G}, \ \ 1 \leq i \leq m,\]
then $\sigma=\tilde{\sigma}$ in $\Omega$. Moreover, there exist a constant $\eta>0$ such that 
for all $\tilde{\sigma}\in W^{1,\infty}(\Omega)$ with 
\begin{equation}\label{smallnesscond}
C_\Omega \parallel \tilde{\sigma} \parallel_{W^{1,\infty}(\Omega)}<\eta,
\end{equation}
where $C_\Omega$ is the best constant in the classical Poincar\'{e} inequality on $H^1_0(\Omega)$. Then the following stability estimate holds 
\begin{equation}
\parallel \sigma -\tilde{\sigma} \parallel_{H^1_0(\Omega)} \leq C \sum_{i=1}^m  \parallel \chi_i \Lambda^*_i (\sigma )- \chi_i \Lambda^*_i(\tilde{\sigma})\parallel_{H^1([0,T] \times \partial B_\rho)},
\end{equation}
where $C>0$ is independent of $\tilde{\sigma}$ and $\sigma$. 
\end{theorem}

\begin{remark}
Notice that $C_\Omega$ is small for a small region $\Omega$ ($\Omega \subset B_r$ for some small $r$), and therefore the condition (\ref{smallnesscond}) is satisfies if the support of $\sigma$ is small in $B_\rho$. This is consistent with the experiments in \cite{BSFECO-JBO-09} where the method is applied on small animals.

\end{remark}

\section{Stability of the Linearized Problem}
\label{sec:lin}
In this section we study the linearized problem associated with (\ref{pde1})-(\ref{measurements}). Fix $\sigma\in C^{\infty}_0(\Omega)$ and let $\tilde{\sigma} \in C^{\infty}_0(\Omega)$. Assume $u_{i}, \tilde{u}_i$ be the corresponding solutions of $(\ref{pde1})$. Then 

\begin{eqnarray}\label{linearized.pde}
\left\{ \begin{array}{ll}
-\Delta \delta u_i+\tilde{\sigma} \delta u_i=-u_i\delta \sigma &\text{in } B_\rho\\
\delta u_i=0 &\text{on }\partial B_\rho, 
\end{array} \right.
\end{eqnarray}
where $\delta u_i= u_i -\tilde{u}_i$ and $\delta \sigma= \sigma- \tilde{\sigma}$. Thus   
\begin{equation*}
\delta u_i(x)=(2\pi)^{-n} \int \int e^{i(x-y)\xi} q(x,\xi) u_i(y) \delta \sigma(y) dy d\xi,
\end{equation*}
where $q(x,\xi)=\frac{-1}{1+\xi ^2}$ mod $\in S^{-3}(B_\rho \times \R^n)$, i.e. 
\[q(x,\xi)+\frac{1}{1+\xi ^2} \in S^{-3}(B_\rho \times \R^n). \]
Recall that a symbol $P(x,\xi)\in S^{m}$ if 
\[|D^{\beta}_x D^{\alpha}_\xi P(x,\xi)| \leq C_{\alpha, \beta} (\sqrt {1+|\xi|^2})^{m-|\alpha|},\]
for all $\alpha, \beta \in \N^n$. See Chapter 1 in \cite{Grigis} for more details about symbols and oscillatory integrals. Now let $\delta p_i(t,x)$ be the solution of the wave equation 

\begin{eqnarray}
\left\{ \begin{array}{ll}
(\partial^2_t-c^2\Delta) \delta p_{i}(t,x)=0 &\text{in } (0,T) \times \R^n\\
\delta p_{i}|_{t=0}=\sigma u_{i}-\tilde{\sigma} \tilde{u}_i \\
\partial_t \delta p_{i}(t,x)|_{t=0}=0.
\end{array} \right.
\end{eqnarray}
Notice that 
\[\sigma u_{i}-\tilde{\sigma} \tilde{u}_i= \Delta (\delta u_i).\]
Hence modulo smooth terms 
\begin{eqnarray*}
\delta p_{i}(t,x)&:=&(2\pi)^{-n}  \sum_{\tau=\pm} \int \int e^{i \varphi_{\tau}(t,x,\xi)} a_{\tau}(t, x,\xi) [\sigma u_{i}(y)-\tilde{\sigma} \tilde{u}_i(y)]d \xi dy\\
&=&(2\pi)^{-n} \sum_{\tau=\pm}  \int \int e^{i \varphi_{\tau}(t,x,\xi)} a_{\tau}(t, x,\xi) \Delta (\delta u_i (y))d \xi dy
\end{eqnarray*}
where the phase function $\varphi_{\pm}$ are homogeneous of order 1 and solve the eikonal equations 
\[\mp \partial_t \varphi_{\pm}=|D_x \varphi_{\pm}|, \ \ \varphi_{\pm}|_{t=0}=x\cdot\xi,\]
 $a_{\pm}$ are amplitudes of order zero satisfying the corresponding transport equations (see equation V.1.50 in \cite{Traves2}). Now define $F_{\pm}$ to be the Fourier integral operators
\begin{equation*}
F_{\pm}(w):=(2\pi)^{-n}  \int \int e^{i \varphi_{\pm}(t,x,\xi)} a_{\pm}(t,x,\xi) w(y) d \xi dy,
\end{equation*}
and $P_i$ to be the pseudodifferential operator 
\begin{equation*}
P_i(w):=(2\pi)^{-n}\int \int e^{i(x-y)\xi} a(x,y,\xi) w(y)dy d \xi,
\end{equation*}
where $a(x,y,\xi)=u_i(y)$ mod $S^{-1}(B_\rho \times B_\rho \times \R^n)$. Then 
\[\delta p_{i}(t,x)=F_{+}P_i(\delta \sigma)+F_{-}P_i(\delta \sigma).\]
Eliminate the dependence on $y$ in the symbol $a(x,y,\xi)$ to get
\begin{equation*}
P_{i}(w):=(2\pi)^{-n}\int \int e^{i(x-y)\xi} b(x,\xi) w(y)dy d \xi,
\end{equation*}
where $b(x,\xi)=u_i(x)$ mod $S^{-1}(B_\rho \times \R^n)$. The composition of the Fourier integral operator $F_{\pm}$ with the pseudodifferential operator  $P_i$ is a Fourier integral operator with the same phase $\varphi_{\tau}$ and amplitude $b_{\tau}(t,x,y,\xi)=u_i(y) a_{\tau}(t,x,\xi)$ mod $S^{-1}(B_\rho \times \R^n)$ (see \cite{Traves2}).  Therefore
\begin{equation}\label{p_i}
\delta p_i(t,x)=(2 \pi)^{-n}\sum_{\tau=\pm}\int \int e^{i \varphi_{\tau}(t,x,\xi)}b_{\tau}(t,x,\xi) u_i(y)\delta \sigma dy d\xi,
\end{equation}
where $b_{\tau}(t,x,\xi)=a_{\tau}(t,x,\xi)$ mod $S^{-1}(B_\rho \times \R^n)$.

The measurements are modeled by the operator 
\begin{equation*}
\Lambda^*(\sigma):=\sum_{i=1}^{m} \chi_i\Lambda^*_i(\sigma)= \sum_{i=1}^{m}\chi_{i} \Lambda (u_{i} \sigma). 
\end{equation*}
Hence 
\begin{eqnarray*}
\Lambda^*(\delta \sigma)=\Lambda^*(\sigma)-\Lambda^*(\tilde{\sigma})&=&\sum_{i=1}^{m}\chi_{i} \Lambda (u_{i} \sigma-\tilde{u}_i\tilde{\sigma})\\
&=& \sum_{i=1}^{m}\chi_{i} \Lambda (\tilde{\sigma} \delta u_i+u_i\delta \sigma)\\
&=&\sum_{i=1}^{m}\chi_{i} \Lambda (\tilde{\sigma} \delta u_i) +\sum_{i=1}^{m}\chi_{i} \Lambda (u_i\delta \sigma).
\end{eqnarray*}
The above pseudodifferential calculus in (\ref{p_i}) indicates that $H_{\delta \sigma}:=\sum_{i=1}^{m}\chi_{i} \Lambda (u_i\delta \sigma)$ is the higher order term in $\Lambda^*(\sigma -\tilde{\sigma})$, and $L_{\delta \sigma}:=\sum_{i=1}^{m}\chi_{i} \Lambda (\tilde{\sigma} \delta u_i)$ may be controlled by $H_{\delta \sigma}$ for small $\delta \sigma$ (see Section 4). Hence we first study invertibility of the operator 
\begin{equation}
p(\delta \sigma):=\sum_{i=1}^{m}(2 \pi)^{-n}\chi_i\sum_{\tau=\pm}\int \int e^{i \varphi_{\tau}(t,x,\xi)}b_{\tau}(t,x,\xi) u_i(y)\delta \sigma dy d\xi=\sum_{i=1}^{m}\chi_{i} \Lambda (u_{i}\delta \sigma)
\end{equation}
and find an approximate inverse. Let $A$ be the back-propagation operator defined in (\ref{A}) and set
\begin{equation}
\kappa (\delta \sigma):=A \left( \sum_{i=1}^{m}\chi_{i} \Lambda (u_{i}\delta \sigma) \right). \\ 
\end{equation}

\begin{proposition}\label{linearized.theo} The operator $\kappa$ is a zero order pseudo-differential operator in a neighborhood of $\Omega$ with principal symbol 
\begin{equation}\label{prinSymb}
\frac{1}{2}\sum_{i=1}^n [\chi_{i}(\gamma_{x,\xi}(\tau_{+}(x,\xi)))+\chi_{i}(\gamma_{x,\xi}(\tau_{-}(x,\xi)))] u_i (x).
\end{equation}
Consequently if (\ref{Uniqueness.con}) and (\ref{stab.con}) hold, then $\kappa$ is an elliptic Fredholm operator on $H^1_0(\Omega)$, and there exists $C>0$ such that
\begin{equation}\label{maininequality}
\parallel \delta \sigma \parallel_{H^1(\Omega)} \leq C \parallel \sum_{i=1}^{m}\chi_{i} \Lambda (u_{i}\delta \sigma) \parallel _{H^1(\mathcal{G})}. 
\end{equation} 
\end{proposition}
{\bf Proof.} First note that  
\[\kappa(\delta \sigma)=\Sigma_{i=1}^{m} A \left(\chi_i\delta p_i(t,x) \right).\]
It follows from Theorem 3 in \cite{SU09} that the principle symbol of $\kappa$ is given by (\ref{prinSymb}). It follows from maximum principle that $u_i \geq \beta>0$ in $\overline{\Omega}$ for all $1\leq i \leq m$. Since $\mathcal{G}$ satisfies (\ref{stab.con}), 
\[\frac{1}{2}\sum_{i=1}^n [\chi_{i}(\gamma_{x,\xi}(\tau_{+}(x,\xi)))+\chi_{i}(\gamma_{x,\xi}(\tau_{-}(x,\xi)))] u_i (x)\geq \frac{\beta}{2}>0.\]
Thus the operator $\kappa$ is elliptic, and therefore it follows from the mapping properties of the back-propagation operator $A$ (see \cite{LLT} and \cite{SU09}) that 
\[\parallel \delta \sigma \parallel_{H^1_0(\Omega)} \leq C \left( \parallel \sum_{i=1}^{m}\chi_{i} \Lambda (u_{i}\delta \sigma) \parallel _{H^1(\mathcal{G})} +\parallel \delta \sigma\parallel_{L^2(\Omega)}\right).\]
Since (\ref{Uniqueness.con}) holds for some $1 \leq j\leq m$,  by Theorem 2 in \cite{SU09}, the measurement $\chi_j \Lambda (\sigma u_j)$  uniquely determines $\sigma u_j$ in $\Omega$. On the other hand, $u_j$ is the unique solution of  
\[-\Delta u_j=-\sigma u_j \ \ u_j|_{\partial B_\rho}=g_j\geq 0.\]
By strong maximum principle we have $u_j>0$ in $\Omega$. Hence if (\ref{Uniqueness.con}) holds, then $\sigma=\frac{\sigma u_j}{u_j}$ is uniquely determined from a single measurement $\chi_j \Lambda (\sigma u_j)$.  Consequently the operator $\sum_{i=1}^{m}\chi_{i} \Lambda (u_{i}\delta \sigma)$ is also injective. Therefore it follows from Proposition V.3.1 in \cite{Tylor} that the estimate (\ref{maininequality}) holds for some constant $C>0$, possibly different from the above constant. \hfill $\Box$

\section{Stability of the Nonlinear Problem}
\label{sec:nlin}
In this section we present the proof of Theorem \ref{MainTheorem}. 
Let $B_{\rho}$ be the ball of radius $\rho$ in $\R^n$ and assume $\sigma$ and $\tilde{\sigma}$ are essentially bounded in $B_\rho$ and $\Omega:=\hbox{supp}(\sigma -\tilde{\sigma})\Subset B_\rho$. For fixed $0 \leq g\in H^{1/2}(\partial B)$ and for $R_i \in \Theta$ define $u_{i}(x)$ to be the unique solution of   
\begin{eqnarray*}
\left\{ \begin{array}{ll}
-\Delta u_i+\sigma u_i=0 &\text{in } B_\rho\\
u=g_{i}&\text{on }\partial B_\rho, 
\end{array} \right.
\end{eqnarray*}
where $g_i(x) =g(R_i x)$, $R_i \in \Theta$. Similarly let $\tilde{u}_i$ be the unique solution of 
\begin{eqnarray*}
\left\{ \begin{array}{ll}
-\Delta \tilde{u}_i+\tilde{\sigma} \tilde{u}_i=0 &\text{in } B_\rho\\
\tilde{u}_i=g_{i}&\text{on }\partial B_\rho. 
\end{array} \right.
\end{eqnarray*}
Then $\delta u_i=u_i-\tilde{u}_i$ satisfies 

\begin{eqnarray}\label{tildsigmaequation}
\left\{ \begin{array}{ll}
-\Delta \delta u_i+\tilde{\sigma} \delta u_i=-u_i\delta \sigma &\text{in } B_\rho\\
\delta u_i=0&\text{on }\partial B_\rho. 
\end{array} \right.
\end{eqnarray}
{\bf Proof of Theorem \ref{MainTheorem}.} Since the mapping $\Lambda: H^1(\Omega) \rightarrow H^1([0,T]\times \partial B_\rho)$ is bounded (see Remark 5 in \cite{SU09}), it is enough to prove the theorem for $\sigma, \tilde{\sigma} \in C_0^{\infty}(\Omega)$. The general result will follow from a standard density argument. 

Multiply (\ref{tildsigmaequation}) by $\delta u_i$ and integrate by parts and use H\"{o}lder's inequality to get 
\begin{equation*}
\parallel \nabla \delta u_i\parallel^2_{L^2(B_\rho)}+\int_{B_\rho}\tilde{\sigma}(\delta u_i)^2dx\leq \parallel u_i\delta \sigma\parallel_{L^2 (B_\rho)} \parallel \delta u_i\parallel_{L^2 (B_\rho)}. 
\end{equation*}
By Poincar\'{e} inequality, there exists $C_\Omega$ such that 
\[ \parallel \delta u_i \parallel_{L^2 (B_\rho)} \leq C_\rho \parallel \nabla \delta u_i \parallel_{L^2 (B_\rho)},\]
where $C_\rho$ is dependent of $\delta u_i$. Thus we have
\begin{equation}
\parallel \nabla \delta u_i\parallel_{L^2(B_\rho)}\leq C_\rho \parallel u_i\delta \sigma\parallel_{L^2 (B_\rho)}.
\end{equation}
Therefore 
\begin{eqnarray*}
\parallel\tilde{\sigma}\delta u_i \parallel_{H^1_0(\Omega)} &\leq& \parallel \tilde{\sigma} \parallel_{W^{1,\infty}(\Omega)} \parallel \nabla \delta u_i\parallel_{L^2(B_\rho)}\\
&\leq &  C_\rho \parallel \tilde{\sigma} \parallel_{W^{1,\infty}(\Omega)}  \parallel u_i\delta \sigma\parallel_{L^2 (B_\rho)}\\
&\leq & C_\rho C_{\Omega} \parallel \tilde{\sigma} \parallel_{W^{1,\infty}(\Omega)}  \parallel \nabla (u_i\delta \sigma)\parallel_{L^2 (B_\rho)},
\end{eqnarray*}
where $C_\Omega$ is the best constant in the classical Poincar\'{e} inequality on $H^1_0(\Omega)$. Thus we have 
\begin{equation}\label{estimate1}
\parallel\tilde{\sigma}\delta u_i \parallel_{H^1_0(\Omega)} \leq  C_\rho C_{\Omega}\parallel \tilde{\sigma} \parallel_{W^{1,\infty}(\Omega)} \parallel u_i \delta \sigma \parallel_{H^1_0(\Omega)}.
\end{equation}

On the other hand, since the mapping $\Lambda: H^1(\Omega) \rightarrow H^1([0,T]\times \partial B_\rho)$ is bounded,

\[ \parallel  \sum_{i=1}^m \chi_i \Lambda (u_i \delta \sigma) \parallel_{H^1([0,T] \times \partial B_\rho)} \leq \tilde{C} \parallel \delta \sigma \parallel_{H^1_0},\]
for some $\tilde{C}>0$. Hence it follows from Proposition \ref{linearized.theo} and (\ref{estimate1}) that 
\begin{eqnarray*}
\parallel \sum_{i=1}^m  \chi_i \Lambda (\sigma u_i- \tilde{\sigma} \tilde{u}_i)\parallel_{H^1([0,T]) \times \partial B_\rho}&=& \parallel \sum_{i=1}^m \chi_i  \Lambda (\tilde{\sigma} \delta{u}_i)+ \chi_i \Lambda(u_i\delta \sigma))\parallel_{H^1([0,T]) \times \partial B_\rho}\\
&\geq & \parallel \sum_{i=1}^m \chi_i  \Lambda(u_i\delta \sigma) \parallel_{H^1([0,T]) \times \partial B_\rho}-\parallel \sum_{i=1}^m \chi_i \Lambda (\tilde{\sigma} \delta{u}_i)\parallel_{H^1([0,T]) \times \partial B_\rho}\\
&\geq & C \parallel \delta \sigma \parallel_{H^1_0(\Omega)}-\tilde{C}\parallel \sum_{i=1}^m \tilde{\sigma} \delta{u}_i\parallel_{H^1_0(\Omega)}\\ 
&\geq & C \parallel \delta \sigma \parallel_{H^1_0(\Omega)} - \tilde{C} C_\rho C_{\Omega} \parallel \tilde{\sigma} \parallel_{W^{1,\infty}(\Omega)}  \sum_{i=1}^{m} \parallel u_i \delta \sigma\parallel_{H^1_0(\Omega)}\\
&\geq & \left( C  - m M_g \tilde{C} C_\rho C_{\Omega} \parallel \tilde{\sigma} \parallel_{W^{1,\infty}(\Omega)}  \right) \parallel \delta \sigma\parallel_{H^1_0(\Omega)},
\end{eqnarray*}
where $M_g$ is the maximum of $g$ on $\partial \Omega$.  
Therefore there exists $\eta>0$ such that if 
\[C_{\Omega}\parallel \tilde{\sigma} \parallel_{W^{1,\infty}(\Omega)} <\eta\]
then 
\begin{equation}
\parallel \sigma -\tilde{\sigma} \parallel_{H^1_0(\Omega)} \leq C^* \parallel \sum_{i=1}^m  \chi_i \Lambda (\sigma u_i- \tilde{\sigma} \tilde{u}_i)\parallel_{H^1([0,T] \times \partial B_\rho)},
\end{equation}
for some $C^*>0$ independent of $\tilde{\sigma}$ and $\sigma$. \hfill $\Box$

\section*{Acknowledgment} GB's work was supported in part by the NSF grant DMS-1408867. AM's work is supported by a start-up grant from University of California, Riverside. The authors would like to thank the anonymous referees for careful reading of the manuscript and helpful comments.

%\bibliography{../../bibliography}
%\bibliographystyle{siam}

\end{document}